\documentclass{amsart}

\usepackage{amssymb}
\usepackage{mathrsfs}

\newtheorem{thm}{Theorem}
\newtheorem{lem}[thm]{Lemma}
\newtheorem{prop}[thm]{Proposition}

\theoremstyle{definition}
\newtheorem{definition}{Definition}

\numberwithin{equation}{section}

\newcommand{\rank}{\mathrm{rank}}
\newcommand{\Pic}{\mathrm{Pic}}
\newcommand{\Hom}{\mathrm{Hom}}

\newcommand{\Lie}{\mathrm{Lie}}

\newcommand{\rmPhi}{\mathrm{\Phi}}

\newcommand{\Sym}{\mathrm{Sym}}

\newcommand{\C}{\mathbb C}
\newcommand{\Z}{\mathbb Z}
\newcommand{\p}{\mathbb P}
\newcommand{\A}{\mathbb A}
\newcommand{\mL}{\mathscr L}
\newcommand{\OO}{\mathcal O}
\newcommand{\g}{\mathfrak g}

\begin{document}

\title{Simple immersions of wonderful varieties}

\author{Guido Pezzini}
\address{Institut Fourier\\ Universit\'e Joseph Fourier\\ B.P. 74\\ 38402 Saint-Martin d'H\`eres\\ France}
\curraddr{Dipartimento di Matematica\\ Universit\`a La Sapienza\\ P.le Aldo Moro 2\\ 00185 Roma\\ Italy}
\email{pezzini@mat.uniroma1.it}

\thanks{Research supported by European Research Training Network {\em LIEGRITS} (MRTN-CT 2003-505078), in contract with CNRS DR17, No 2. }

\date{June 23, 2006}

\begin{abstract}
Let $G$ be a semisimple connected linear algebraic group over $\C$, and $X$ a wonderful $G$-variety. We study the possibility of realizing $X$ as a closed subvariety of the projective space of a simple $G$-module. We describe the wonderful varieties having this property as well as the linear systems giving rise to such immersions. We also prove that any ample line bundle on a wonderful variety is very ample.
\end{abstract}

\maketitle

\section{Introduction}

Wonderful varieties are projective algebraic varieties (for us, over the field of complex numbers $\C$) endowed with an action of a semisimple connected algebraic group $G$, having certain properties which have been inspired by the compactifications of symmetric homogeneous spaces given by De Concini and Procesi in \cite{DP83}. Wonderful varieties turn out to have a significant role in the theory of spherical varieties, which are a class of $G$-varieties representing a common generalization of flag varieties and toric varieties.

In this paper we answer a question raised by Brion at the end of the paper \cite{Br90}. Given a wonderful $G$-variety $X$, we want to study whether there exists a $G$-equivariant closed immersion in a projective space $\p(V)$ where $V$ is a simple $G$-module (a ``simple immersion'', for brevity). This fact is true if $X$ is a complete symmetric variety in the sense of \cite{DP83}, but is false for other easy examples, such as $\p^1\times\p^1$ under the diagonal action of $PSL_2$.

In our main theorem (theorem \ref{thm:main}) we prove a necessary and sufficient condition for this to be true, and find all the linear systems which give rise to such an immersion. The condition is given in terms of the stabilizers in $G$ of the points of the variety; it can also be stated in terms of some known invariants of wonderful varieties (the {\em spherical roots}).

Our approach consists in reducing the problem to a small family of wonderful varieties: those of rank $1$ which are not a parabolic induction. This family is finite and classified for any given $G$, see the works of Ahiezer (\cite{Ah83}), Huckleberry and Snow (\cite{HS82}), Brion (\cite{Br89a}), so we can carry on the proof case-by-case.

The same technique is used to show that on a wonderful variety any ample line bundle is very ample.

\subsection*{Acknowledgements}
The Author thanks Prof.~M.~Brion for all his precious help in the developing of this work, and Prof.~D.~Luna for the fruitful discussions on the subject.

\section{Definitions}\label{sect:def}
\subsection{Wonderful varieties}
Throughout this paper, $G$ will be a semisimple connected algebraic group over $\C$. We also suppose $G$ simply connected. In $G$ we fix a Borel subgroup $B$, a maximal torus $T\subset B$, we denote by $\rmPhi$ the corresponding root system, by $\rmPhi^+$ the positive roots and by $S$ the set of simple roots. Also, we denote by $B_-$ the Borel subgroup opposite to $B$, i.e.~such that $B\cap B_- = T$. Given a subset $S'$ of the simple roots, we denote by $\rmPhi_{S'}$ the associated root subsystem, and we denote by $G_{S'}$ (resp.~$G_{-S'}$) the associated parabolic subgroup containing $B$ (resp.~$B_-$). We will also use the notation $\C^\times$ to denote the multiplicative group $\C\setminus\{0\}$.

\begin{definition} \cite{Lu01}
A {\em wonderful $G$-variety} is an irreducible algebraic variety $X$ over $\C$ such that:
\begin{enumerate}
\item X is smooth and complete;
\item $G$ has a open (dense) orbit on $X$, and the complement is the union of ($G$-stable) prime divisors $D_i$ ($i = 1,\ldots,r$), which are smooth, with normal crossings and satisfy $\bigcap_{i=1}^r D_i \neq \emptyset$;
\item If $x,x'\in X$ are such that $\left\{ i \;|\; x\in D_i\right\} = \left\{ i \;|\; x'\in D_i\right\}$, then $x$ and $x'$ lie on the same $G$-orbit.
\end{enumerate}
The number $r$ of $G$-stable prime divisors is the {\em rank} of $X$.
\end{definition}

A wonderful variety $X$ is always spherical, i.e.~a Borel subgroup has an open dense orbit on $X$ (see \cite{Lu96}). We can introduce some data associated to $X$ coming from the theory of spherical varieties, and fix some notations (for details, see \cite{Kn96}, \cite{Lu01}):

\begin{enumerate}
\item $H = $ the stabilizer of a point in the open $G$-orbit of $X$, so that this orbit is isomorphic to $G/H$: it is known that $H$ has finite index in its normalizer $N_GH$; we choose the point stabilized by $H$ to be also in the open $B$-orbit, so that $BH$ is open in $G$;
\item $\Xi_X = \left\{B\textrm{-weights of functions in } \C(X) \textrm{ which are }B\textrm{-eigenvectors}\right\}$, where $G$ acts on rational functions on $X$ in the usual way: $(gf)(x)=f(g^{-1}x)$, and the weight of a function $f$ is the character $\chi\colon B\to \C^\times$ such that $bf=\chi(b)f$ for all $b\in B$;
\item $\Delta_X = \left\{ B\textrm{-stable but not } G\textrm{-stable prime divisors of }X \right\}$, whose elements are called the {\em colours} of $X$;
\item for any colour $D$, we define $\rho_X(D) \in \Hom_\Z(\Xi_X,\Z)$ in the following way: $\langle\rho_X(D),\chi\rangle=\nu_D(f_\chi)$ where $\nu_D$ is the discrete valuation on $\C(X)$ associated to $D$, and $f_\chi$ is a rational function on $X$ being a $B$-eigenvector with weight $\chi$; this functional $\rho_X(D)$ is well defined because $X$ has an open $B$-orbit and thus any weight $\chi$ determines $f_\chi$ up to a multiplicative constant;
\item for any simple root $\alpha\in S$, we say that $\alpha$ ``moves'' a colour $D\in\Delta_X$ if $D$ is non-stable under the action of $G_{\{\alpha\}}$;
\item $z = $ the unique point fixed by $B_-$ on $X$; it lies on $Z=Gz$ the unique closed $G$-orbit;
\item $\Sigma_X=\left\{ T\textrm{-weights of the } T\textrm{-module } T_zX/T_z(Gz)\right\}$; the elements of $\Sigma_X$ are called the {\em spherical roots} of X, and the cardinality of $\Sigma_X$ is equal to $\rank X$;
\item $P_X =$ the stabilizer in $G$ of the open $B$-orbit; this is a parabolic subgroup containing $B$;
\item $S^p_X=$ the subset of simple roots associated to the parabolic subgroup $P_X$, so that in our notations: $P_X = G_{S^p_X}$.
\end{enumerate}

\subsection{Parabolic induction} \label{subsect:induction}
Let $X$ be a wonderful $G$-variety, and suppose that the stabilizer $H$ of a point in the open $G$-orbit is such that $R(Q) \subseteq H \subseteq Q$ for some parabolic subgroup $Q$ of $G$, where $R(Q)$ is the radical of $Q$. Then $X$ is isomorphic to $G\times_Q Y$ where $Y$ is a $Q$-variety where the radical $R(Q)$ acts trivially. Moreover, $Y$ turns out to be wonderful under the action of $Q/R(Q)$, thus also under the action of $L$ a Levi subgroup of $Q$. Here $G\times_Q Y$ is defined as the quotient $G\times Y/\sim$ where $(g,x)\sim(gq,q^{-1}x)$ for all $q\in Q$.

\begin{definition}
Such a wonderful variety $X\cong G\times_Q Y$ is said to be a {\em parabolic induction} of $Y$ by means of $Q$. A wonderful variety which is not a parabolic induction is said to be {\em cuspidal}.
\end{definition}

We will need some facts upon wonderful varieties which are parabolic induction. The main idea is that the $L$-action on $Y$ determines the whole structure of $X$, and this is the reason why the study of wonderful varieties most often reduces to the study of cuspidal ones.

It is convenient to make some choice upon $Q$ and $L$, using conjugation whenever necessary: we choose $Q$ so that it contains $B_-$, and let $S(Q)\subset S$ be the subset of simple roots associated to $Q$. We can choose $L$ such that $B\cap L$ is a Borel subgroup of $L$.

Let $\phi\colon X\to G/Q$ be the map given by $[g,y]\mapsto gQ$; we can identify $Y$ with $\phi^{-1}(Q)$. There are some colours of $X$ that go surjectively on $G/Q$ via $\phi$: such a colour is equal to $\overline{B D}$ for a colour $D$ of $Y$. The other colours of $X$ are the pull-back of the colours of $G/Q$ along $\phi$. Therefore $\Delta_X$ can be identified with the disjoint union of $\Delta_{G/Q}$ and $\Delta_{Y}$; however, in order to avoid confusion, we use the notation $\widetilde D$ to denote the element in $\Delta_X$ associated to $D\in\Delta_{G/Q}$ or $D\in\Delta_Y$.

With this identification, any simple root moving some colour coming from $Y$ must belong to $S(Q)$, while simple roots moving colours coming from $G/Q$ must be in $S\setminus S(Q)$.

\subsection{Line bundles}\label{subsect:linebundles}
Since $X$ is spherical and has only one closed $G$-orbit, we have the following description of $\Pic(X)$:

\begin{prop}\cite{Br89}\label{prop:bunsphe}
The Picard group of $X$ has a basis consisting of the classes of the colours. Moreover, the divisors which are generated by global sections (resp.~ample) are the linear combinations of these classes having non-negative (resp.~positive) coefficients.
\end{prop}

Any line bundle $\mL$ on $X$ has a unique $G$-linearization, that is, a $G$-action on the total space of the bundle such that the projection on $X$ a $G$-equivariant map, and such that it is a linear action on the fibers (see \cite{KKLV89}). This determines an action of $G$ on the space of global sections $\Gamma(X,\mL)$. Our $X$ is spherical, so this $G$-module has no multiplicities, which means that any simple $G$-module appears no more than once (see \cite{Br97}). Moreover, the highest weights of the simple $G$-modules which actually appear can be described precisely.

Let $\mL$ be a line bundle on $X$, and suppose it is associated to an effective divisor of the form $\delta=\sum_{D\in\Delta_X}n_D D$ with $n_D\geq 0$ for all $D$. We will also use the standard notation $\OO(\delta)$ for such a line bundle. Consider its canonical section $\sigma_\mL\in\Gamma(X,\mL)$; since the colours are $B$-stable then $\sigma_\mL$ is $B$-proper: call its $B$-weight $\chi_\mL$. Then, the highest weights of the simple modules appearing in $\Gamma(X,\mL)$ are the dominant weights which can be expressed as $\chi_\mL + \xi$, where $\xi$ is a linear combination of spherical roots with non-positive coefficients, and $\langle\rho_X(D),\xi\rangle + n_D \geq 0$ for all colours $D$ of $X$.

This result is established in \cite{Br89}, after some analysis of the action of $G$ on $\Gamma(X,\mL)$ in relation with the usual induced action on $\C(X)$. It is useful to recall here the main idea; we start from the following equality of vector spaces:
\[ \Gamma(X,\mL) = \left\{ f\in\C(X)\;\;|\;\;(f)+\delta\geq 0\right\}. \]
The restriction to the open $G$-orbit $G/H\subseteq X$ induces an inclusion:
\[ \Gamma(X,\mL) \subseteq \Gamma(G/H,\mL). \]
The sections in $\Gamma(G/H,\mL)$ are quotients $f_1/f_2$ of regular functions on $G$ which are $H$-proper (with same weight) under right translation, and such that the zeros of $f_2$ are ``less than or equal to'' $\delta$. Let $f_\delta\in\C[G]$ be a global equation of $\delta$ pulled back on $G$ via the projection $G\to G/H$ ($f_\delta$ is unique up to a multiplicative constant). It is a regular function on $G$ and it is $H$-proper under right translation; denote $\lambda$ its $H$-weight and $\C[G]^{(H)}_{\lambda}$ the set of all $H$-proper functions with that weight. We have:
\[ \Gamma(G/H,\mL) = \left\{ \frac{f}{f_\delta} \;\;|\;\; f\in\C[G]^{(H)}_{\lambda} \right\};\]
and the map $f/f_\delta\mapsto f$ gives a inclusion $\Gamma(X,\mL)\subseteq \C[G]^{(H)}_{\lambda}$, where the latter in turn provides the $G$-module structure of $\Gamma(X,\mL)$ via the left translation action of $G$. In this way, the canonical section $\sigma_\mL$ is represented by the constant rational function $1$, and it corresponds to the highest weight vector $f_\delta$.

Notice that if $g\in G$ fixes $f_\delta$ under left translation, then $g$ acts on a section of $\mL$ in the same way as it acts on the corresponding rational function with the usual action induced on $\C(X)$.

\section{Simple immersions}
\subsection{Main theorem}
The basic examples are the two wonderful $SL_2$-varieties: $\p^1\times\p^1$ and $\p^2\cong\p(\Sym^2(\C^2))$. The former does not admit any immersion into the projective space of a simple $SL_2$-module, whereas $\p^2$ does (and even happens to be such a projective space), as we will see in \ref{ssect:nonrigid}. Notice that the open orbits of these two varieties are isomorphic resp.\ to $SL_2/T$ and $SL_2/N_{SL_2}T$, with $|N_{SL_2}T/T|=2$. 

\begin{thm}
\label{thm:main}
Let $X$ be a wonderful $G$-variety. There exist a simple $G$-module $V$ and a $G$-equivariant closed immersion $X\to\p(V)$ if and only if the stabilizer of any point of $X$ is equal to its normalizer. For any fixed $V$, this immersion is unique.
\end{thm}

Varieties satisfying the condition of this theorem are also called {\em strict}. The proof of the theorem will take place in section \ref{sect:proof}. We will also see in that section that we can characterise all simple modules admitting such an immersion for a given $X$.

It is evident that any map as in the theorem is given by a linear system which corresponds to a simple submodule of $\Gamma(X,\mL)$ for some ample line bundle $\mL$.

The easiest case is where $X$ has rank zero. Indeed, rank zero wonderful varieties are exactly the generalized flag varieties $G/Q$, $Q$ a parabolic subgroup of $G$. A part of the classical Borel-Weil theorem states that the space of global sections of any ample line bundle on $G/Q$ is an irreducible $G$-module, and it gives a closed immersion in the corresponding projective space.

For $X$ of any rank, $\Gamma(X,\mL)$ is not irreducible; however we have strong restrictions on the submodules we can consider:

\begin{prop}\label{prop:notall}
Let $\mL$ be generated by global sections, and $V$ be a simple submodule of $\Gamma(X,\mL)$. If the associated rational map $F\colon X\dashrightarrow\p(V^*)$ is regular and birational onto its image, then the centre of $G$ acts trivially on $X$ and the highest weight of $V$ is the ``highest possible'' among the simple modules occurring in $\Gamma(X,\mL)$, which is the weight $\chi_\mL$.
\end{prop}
\begin{proof}
The assertion on the centre of $G$ is evident, since an element in the centre acts as a scalar on $V$ and thus trivially on $\p(V)$.

Let $Z$ be the unique closed $G$-orbit on $X$: the restriction of $\mL$ to $Z$ gives a $G$-equivariant linear map $\Gamma(X,\mL)\to\Gamma(Z,\mL|_Z)$, and $\Gamma(Z,\mL|_Z)$ is a simple $G$-module by the Borel-Weil theorem. Take a point $z\in Z$ fixed by $T$, such that its stabilizer $Q$ is a parabolic subgroup with $BQ$ open in $G$. Then, $Q$ acts on the fiber at $z$ of $\mL$ with the character $-\chi_\mL$, thus the highest weight of $\Gamma(Z,\mL|_Z)$ is $\chi_\mL$, again by the Borel-Weil theorem.

Therefore the map $\Gamma(X,\mL)\to\Gamma(Z,\mL|_Z)$ is zero on all simple submodules of $\Gamma(X,\mL)$, except for the one having weight $\chi_\mL$. Recall that all simple submodules appearing have highest weight of the form $\chi_\mL + \xi$, where $\xi$ is a linear combination of spherical roots with nonpositive coefficients: we have shown that all simple submodules with $\xi\neq 0$ cannot give a regular map, since their sections restrict to zero on $Z$.
\end{proof}

\begin{definition}
We denote $V_\mL$ the simple submodule of $\Gamma(X, \mL)$ having highest weight $\chi_\mL$.
\end{definition}

The following proposition is related to our problem, and will be useful.

\begin{prop}\label{prop:brionemb} \cite{Br97}
Let $X$ be a wonderful variety: there exist a simple $G$-module $M$ and a vector $v\in M$ such that:
\[ H \subseteq G_{[v]} \subseteq N_GH \]
where $[v]\in\p(M)$ is the point corresponding to $v$ and $G_{[v]}$ is its stabilizer in $G$. Our $X$ is isomorphic to the normalization of the variety $\overline{G[v]}\subseteq \p(M)$ in the field $\C(G/H)$, which contains $\C(G/G_{[v]})$.
\end{prop}

Notice that if $\mL$ is generated by its global sections, then the global sections belonging to $V_\mL$ suffice to generate $\mL$, as one can see from the proof of proposition \ref{prop:notall}.

It is useful to state a slightly different version of theorem \ref{thm:main}, taking also into account proposition \ref{prop:notall}.

\begin{thm}\label{thm:SI}
Let $X$ be a wonderful $\overline{G}$-variety, where $\overline{G}$ is the adjoint group of $G$. Let $\mL$ be an ample line bundle on $X$. Then the simple $G$-submodule $V_\mL$ of $\Gamma(X,\mL)$ gives a closed immersion $F_\mL\colon X\to \p\left(V_\mL^*\right)$ if and only if the following condition holds:
\begin{itemize}
\item[(R)]any wonderful $G$-subvariety $X'$ of rank $1$ of $X$ is rigid, i.e.\ its generic stabilizer $H'$ is equal to its normalizer $N_G(H')$.
\end{itemize}
which is also equivalent to the following combinatorial condition:
\begin{itemize}
\item[(R')]for any spherical root $\gamma$ of $X$, there exist no rank $1$ wonderful $\overline{G}$-variety $X'$ having spherical root $2\gamma$ and such that $S^p_X=S^p_{X'}$.
\end{itemize}
\end{thm}

The theorem follows from lemma \ref{lemma:reduction} and from section \ref{sect:rank1}. The equivalence of the two conditions (R) and (R') is easy: the spherical roots of $X$ and its wonderful subvarieties of rank $1$ are in bijection, in such a way that the subvariety $X^\gamma$ associated to $\gamma\in\Sigma_X$ satisfies $\Sigma_{X^\gamma}=\{\gamma\}$, and $S^p_{X^\gamma}=S^p_X$ (see \cite{Lu01}). Now the equivalence follows at once from the classification of rank $1$ varieties.

\subsection{Reduction to rank $1$}
Let $X$ be a wonderful $\overline G$-variety of any rank, with an ample line bundle $\mL = \OO(\delta)$ where $\delta = \sum_{D\in\Delta_X} n_D D$ ($n_D >0$ for all $D$).

In the case of $X$ being of rank zero our main result is immediate, as we have already noticed, thanks to the Borel-Weil theorem.

The general case can be reduced to the study of the rank one case. We begin with the following lemma, which can be essentially found in \cite{Lu02}:
\begin{lem} \cite{Lu02} \label{lemma:luna}
The following conditions are equivalent:
\begin{enumerate}
\item $F_\mL$ is a closed immersion;
\item the application $T_z F_\mL$ between tangent spaces induced by $F_\mL$ in a point $z$ of the closed $G$-orbit $Z$ of $X$ is injective;
\item the restrictions of $F_\mL$ to the rank $1$ wonderful sub-$G$-varieties of $X$ are closed immersions.
\end{enumerate}
\end{lem}
The rank $1$ wonderful sub-$G$-varieties of $X$ are the intersections of any $r-1$ prime divisors stable under $G$ (the $D_i$'s of the definition) where $r = \rank X$.
\begin{proof}
The implications (1)$\Rightarrow$(2) and (1)$\Rightarrow$(3) are obvious. We begin with (2)$\Rightarrow$(1).

The map $F_\mL$ is $G$-equivariant and $Z$ is the unique closed $G$-orbit, therefore (2) ensures that $T_x F_\mL$ is injective for all $x\in X$. This implies also that $F_\mL$ is finite, a consequence of the Stein factorization and the finiteness of the fibers of $F_\mL$.

Let $Z'=F_\mL(Z)$: since $F_\mL$ is finite, $Z\to Z'$ is an isomorphism. We also have that $F_\mL^{-1}(Z')=Z$ since $Z$ is the unique closed $G$-orbit. The set $\{x\in X \;\;|\;\; F_\mL^{-1}(F_\mL(x))=\{x\}\}$ is a $G$-stable open subset of $X$, and it contains $Z$: therefore it is equal to the whole $X$; this shows that $F_\mL$ is a closed immersion.

We show that (3)$\Rightarrow$(2). For all spherical roots $\gamma$, let $(T_z X)_\gamma$ be the subspace of $T_z X$ where $T$ acts with weight $\gamma$. The tangent space $T_z X$ is the sum of the $(T_z X)_\gamma$'s (which are in direct sum), plus $T_z Z$ (with possibly non-trivial intersection with the previous subspaces).

For all spherical roots $\gamma$ we have that $(T_z X)_\gamma + T_z Z$ is the tangent space at $z$ of some rank $1$ wonderful subvariety $X^\gamma$ of $X$, and (3) ensures that $T_z F_\mL |_{T_z X^\gamma}$ is injective for all $\gamma$. This implies (2).
\end{proof}

Our approach will derive from the following lemma. Here we fix a representative $\dot{w}_0\in N_GT$ of the longest element in the Weyl group $N_GT/T$; recall that $P_X$ is the stabilizer of the open $B$-orbit of $X$.

\begin{lem}\label{lemma:diffcond}
Let $d$ be the dimension of $X$, and let $\sigma_\mL$ be the canonical section of $\mL$, viewed as the constant rational function $1$ on $X$. Then the simple submodule $V_\mL$ of $\Gamma(X,\mL)$ gives a closed immersion $F_\mL\colon X\to \p\left(V_\mL^*\right)$ if and only if there exist $u_1,\ldots,u_d \in R^u(P_X)$ such that the Jacobian matrix of the functions%
\footnote{The Reader should make no confusion here: we consider the sections of $\mL$ as rational functions on $X$, but $G$ does not act on $\Gamma(X,\mL)$ via the usual action on $\C(X)$ (see the end of section \ref{sect:def}).}
$(u_1\dot{w}_0)\sigma_\mL,\ldots, (u_d\dot{w}_0)\sigma_\mL$ is nondegenerate at $z$.
\end{lem}
\begin{proof}
Since $V_\mL$ is simple then $G\sigma_\mL$ spans the whole $V_\mL$. This is true also if we replace $G$ with any non-empty open set, like $R^u(P_X) \dot{w}_0 P_X$. The section $\sigma_\mL$ is an eigenvector under the action of $P_X$, so we have that $(R^u(P_X) \dot{w}_0) \sigma_\mL$ spans $V_\mL$.

Therefore there is no harm in supposing that the map $F_\mL$ is given (in coordinates) by sections of the form $(u \dot{w}_0) \sigma_\mL$, for $u\in R^u(P_X)$. Now it follows from lemma \ref{lemma:luna} that the problem of having a closed immersion is local in $z$ and can be checked just on the induced application on the tangent spaces.
\end{proof}

The advantage of this pont of view is that we can use the {\em canonical chart} $X_{Z,B}$, which is the open set of $X$ where $\sigma_\mL$ is non-zero. The notation $X_{Z,B}$ comes from the following equivalent definitions (see \cite{Br97}):
\[ X_{Z,B} = \left\{ x\in X\;\;|\;\; \overline{Bx}\supseteq Z \right\} = X\setminus\bigcup_{D\in\Delta_X}D\]
This is an affine open set; its stabilizer in $G$ is exactly $P_X$. The canonical chart has a very useful description:

\begin{prop} \label{prop:canonicalchart} \cite{Br97}
Let $X$ be any wonderful $G$-variety, and let $L$ be a Levi subgroup of $P_X$. We can suppose that $L$ contains $T$. There exists a $L$-stable closed affine subvariety $M\subseteq X_{Z,B}$ which intersects $Z$ exactly in $z$, and such that we have a $P_X$-equivariant isomorphism:
\[
\begin{array}{ccc}
R^u(P_X)\times M & \longrightarrow & X_{Z,B} \\
  (u,x)     & \longmapsto     & ux
\end{array}
\]
where $P_X$ acts on the product in the following way: if we write an element in $P_X$ as $vl$ where $v\in R^u(P_X)$ and $l\in L$, then $vl(u,x)=(vlul^{-1},lx)$. Moreover, $M$ is an affine space where $(L,L)$ acts trivially and the whole $L$ acts linearly with weights the spherical roots of $X$.
\end{prop}

We are now able to reduce the problem to rank $1$ varieties.

\begin{lem}\label{lemma:reduction}
If theorem \ref{thm:SI} holds for all cuspidal wonderful varieties of rank $1$, then it holds for all wonderful varieties.
\end{lem}
\begin{proof}
We maintain the notations of the theorem. We first prove that the theorem can be reduced to wonderful varieties of rank $1$, and then we reduce to the cuspidal case. Hence now we suppose that the theorem holds for all rank $1$ wonderful varieties.

Let $X$ be a wonderful variety of rank $r$ and suppose that it satisfies the condition (R) of the theorem. 

A spherical root $\gamma\in\Sigma_X$ is associated to a unique $G$-stable wonderful subvariety $X^\gamma\subseteq X$ of rank $1$, having spherical root $\gamma$ and such that $S^p_{X^\gamma}=S^p_X$: the condition (R) of the theorem holds for all these $X^\gamma$ too.

Let $\mL$ be an ample line bundle on $X$, and consider $\mL^\gamma:= \mL|_{X^\gamma}$. The hypotesis of our lemma applies, and the map $F_{\mL^\gamma}\colon X^\gamma\to \p(V_{\mL^\gamma}^*)$ associated to $\mL^\gamma$ is a closed immersion. Thanks to lemma \ref{lemma:diffcond}, this fact can be expressed in terms of the Jacobian matrix at $z$ of some $G$-translates of the canonical section of $\mL^\gamma$. This section is actually $\sigma_\mL$ restricted to $X^\gamma$, so we have proven that $(F_\mL)|_{X^\gamma}$ is a closed immersion of $X^\gamma$ into $\p(V_\mL^*)$, for all $\gamma$.

We conclude that $F_\mL$ is a closed immersion of the whole $X$ thanks to lemma \ref{lemma:luna}.

On the contrary, suppose that $X$ does not satisfy the condition (R) of the theorem, for some $\gamma\in\Sigma_X$. Consider $X^\gamma$ as before: it does not satisfy the condition (R) and so it can't be embedded in the projective space of any simple $G$-module thanks to the hypothesis of this lemma. Thus $\mL$ here cannot give a closed immersion $X\to\p(V_\mL^*)$ because it would restrict to a closed immersion of $X^\gamma$.

We have reduced the problem to rank $1$ varieties, it remains to reduce the problem to the cuspidal case.

Suppose that the theorem holds for all cuspidal wonderful varieties of rank $1$, and let $X$ be a non-cuspidal wonderful variety of rank $1$. Thus $X=G\times_Q Y$ where $Q$ is a proper parabolic subgroup of $G$ and $Y$ is a cuspidal wonderful $L$-variety of rank $1$ for $L$ a Levi subgroup of $Q$. We choose $Q$ and $L$ as in \ref{subsect:induction}, and take an ample line bundle $\mL$ on $X$.

Define\footnote{If we have a group $\Gamma$ acting on a set $A$ then we use the standard notation $A^\Gamma$ to denote the fixed points of $\Gamma$ in $A$.} $W= V_\mL^{R^u(Q)}$, and observe that $Y=X^{R^u(Q)}$. The map $F_\mL$ sends $Y$ into $\p(W^*)$, and $W$ is a simple $L$-module.

Also, $X$ satisfies the condition (R) of the theorem if and only if $Y$ does. If $F_\mL$ is a closed immersion of $X$ into $\p(V_\mL^*)$ then its restriction to $Y$ is a closed immersion of $Y$ into $\p(W^*)$, associated to the ample line bundle $\mL|_Y$.

Viceversa, suppose that $\mL|_Y$ is a closed immersion of $Y$ into $\p(W^*)$; this implies that $T_z F_\mL$ is injective if restricted to $T_z Y$. We also know thanks to the Borel-Weil theorem that $T_z F_\mL$ is injective if restricted to $T_z Z$.

The tangent space of $X$ at $z\in Z$ can be written as: $T_zX= T_zY + T_z Z$ (with non-zero intersection); we want to conclude that $T_z F_\mL$ is injective on all $T_zX$. We use that $X= G\times_Q Y$: the quotient $T_zX/T_zY$ is $T$-isomorphic to $T_P(G/Q)$; the $T$-weights appearing in it do not belong to the span of the roots of the Levi factor $L\subset Q$. On the other hand, the $T$-weights appearing in $T_zY$ are those corresponding to the tangent space of the closed $L$-orbit of $Y$, with in addition the spherical root of $Y$: all of them lie in the weight lattice of $L$. Therefore none of the $T$-weights appearing in $T_zX/T_zY$ can appear also in $T_zY$; since $T_z F_\mL$ is $T$-equivariant, we conclude that it is injective on all $T_zX$.

Lemma \ref{lemma:luna} applies again and we conclude that $F_\mL$ is a closed immersion of $X$ into $\p(V_\mL^*)$.
\end{proof}

\section{Rank $1$ cuspidal wonderful varieties}\label{sect:rank1}
\subsection{General considerations}\label{subsect:general}
Now let $X$ be a rank $1$ wonderful $\overline G$-variety, having dimension $d$ and spherical root $\gamma$, with the ample line bundle $\mL$ as in the previous section.

The idea here is to use lemma \ref{lemma:diffcond} on the canonical chart $X_{Z,B}$.

Proposition \ref{prop:canonicalchart} describes $X_{Z,B}$ as the product of $Z\cap X_{Z,B} = R^u(P_X)$ and an affine space $M$, of dimension equal to $\rank X=1$. The maximal torus $T$ acts linearly on $M$ via the spherical root of $X$.

Therefore $X_{Z,B}\cong \A^{\dim R^u(P_X)+1}=\A^d$. Our strategy is based on the fact that if we can find the section $\dot{w}_0\sigma_\mL$ explicitly as a regular function on $X_{Z,B}$, then the condition of lemma \ref{lemma:diffcond} can be examined when working only on $X_{Z,B}$.

Here it is convenient to consider the action of $P_X$ on $\C[X_{Z,B}]\supseteq\Gamma(X, \mL)$ induced by the one on $\C(X)$, instead of that induced by the linearization of $\mL$. This causes no harm, since $R^u(P_X)$ and $(L,L)$ fix the canonical section $\sigma_\mL$, thus the two actions are the same for $R^u(P_X) (L,L)$. The difference in the action of $T$ is simply the shift by $\chi_\mL$ of all weights (see the end of section \ref{sect:def}).

Fix a global coordinate system $x_1, \ldots, x_{d-1}, y$ on $X_{Z,B}$ with the $x_i$'s relative to $R^u(P_X)$ and the $y$ relative to $M\cong \A^1$. Choose the $x_i$'s to be $T$-eigenvectors; this can be accomplished using the decomposition in root spaces $\mathfrak g_\alpha \subseteq \Lie(R^u(P_X))$  and the isomorphism $\exp\colon \Lie(R^u(P_X)) \to R^u(P_X)$. We will use the corresponding subset of positive roots as an alternative set of indexes for our variables: 
\[\{x_1, x_2, \ldots, x_{d-1} \} = \{ x_\alpha \;|\; \g_\alpha\subseteq \Lie(R^u(P_X))\} = \{ x_\alpha \;|\; \alpha \in \rmPhi^+\setminus\rmPhi_{S^p_X} \} .\]
Notice that the weight of the function $x_\alpha$ is the negative root $-\alpha$.

If we think of $\dot{w}_0\sigma_\mL$ as a regular function on $X_{Z,B}$, then it is a polynomial in these coordinates. Its zero locus is $\dot{w}_0\delta = \sum_{D\in\Delta_X}n_D (\dot{w}_0 D)$ intersected with $X_{Z,B}$.

We have:
\[
\dot{w}_0 \sigma_\mL = \prod_{D\in\Delta_X} \left(\dot{w}_0 \sigma_{\OO(D)}\right)^{n_D}
\]
It is convenient to write it also as a polynomial in $y$:
\begin{equation} \label{formula:w0sigmaL} \dot{w}_0\sigma_\mL = f_m(x) y^m + \ldots + f_1(x) y + f_0(x) \end{equation} 
where $x = (x_1,\ldots,x_{d-1})$. The function $y$ is the equation of $R^u(P_X)$ inside $X_{Z,B}$. So in particular $y$ is a $T$-eigenvector, with weight $-\gamma$.

Since $\dot{w}_0 \sigma_\mL$ is a $T$-eigenvector too, each summand $f_i(x)y^i$ must be a $T$-eigenvector with the same weight. This implies that $f_i(x)$ is a $T$-eigenvector, with weight a sum of negative roots.

Finally, $f_0(x)$ is the equation (in $X_{Z,B}$) of $\dot{w}_0\delta|_{Z}$. The latter is the sum (with multiplicities) of some of the $B_-$-stable prime divisors of $Z$. The explicit equations of these prime divisors can be easily found; anyway often the following lemma suffices:

\begin{lem}\label{lemma:weightsonZ}
If $E$ is a colour of $Z$ and it is moved by $\alpha_i\in S$, then the $T$-weight of the equation of $\dot{w}_0 E$ on $R^u(P_X)$ is $\dot{w}_0(\omega_i)-\omega_i$ where $\omega_i$ is the fundamental dominant weight associated to $\alpha_i$.
\end{lem}
\begin{proof}
Let $P_-=G_{-S^p_X}$ be the opposite subgroup of $P_X$ with respect to $T$: it is also the stabilizer of $z$ the unique point fixed by $B_-$ (see for example \cite{Br97}). Consider the pullback along $\pi:G\to G/P_- = Z$ of the colour $E$. Call $f_E$ its global equation on $G$: it is a $B$-eigenvector under the left translation action, of weight $\omega_i$. 

The map $\pi$ sends $R^u(P_X)$ isomorphically onto the the canonical chart of $Z$. The equation of $\dot{w}_0 E$ on $R^u(P_X)$ corresponds then to the rational function $\frac{\dot{w}_0 f_E}{f_E}$ on $Z$ (where $\dot{w}_0$ acts on regular functions by left translation), and the lemma follows.
\end{proof}

The multiplicities of the intersections between colours and $Z$ are given by the following proposition due to D.~Luna:

\begin{prop} \label{prop:mult}
Let $X$ be a wonderful variety (of any rank), with closed orbit $Z$ and let $D$ be a colour. Then $D$ intersects $Z$ with multiplicity $1$, unless $D$ is moved by $\alpha\in S$ such that $2\alpha\in\Sigma_X$; in this case $D$ intersects $Z$ with multiplicity $2$.
\end{prop}
\begin{proof}
We use the Poincar\'e duality between divisors and curves. Using a Bialynicki-Birula decomposition of $X$ (see \cite{BB73}, \cite{BB76}), one can find a set of $B_-$-stable curves which are a dual basis of the basis of $\Pic(X)$ given by the colours. That is, for each $D\in\Delta_X$ there exists a $B_-$-stable curve $c_D$ in $X$ such that:
\begin{itemize}
\item[--] $c_D\cap D' = \emptyset$ if $D\neq D'$,
\item[--] $c_D$ and $D$ intersect (transversally) in only one point.
\end{itemize}
On the other hand, in \cite{Br93} it is shown that all $B_-$-stable curves are the following:
\begin{enumerate}
\item curves $c_\alpha$ for $\alpha\in S\setminus S^p_X$, which are in $Z$ and are dual to the colours of $Z$;
\item curves $c_\alpha^+$ and $c_\alpha^-$, for $\alpha\in S\cap\Sigma_X$;
\item curves $c_\alpha'$ for $\alpha\in S\cap\frac{1}{2}\Sigma_X$.
\end{enumerate}
See \cite{Lu97} for a description of these curves. In the Chow ring (or, equivalently, in the cohomology) of $X$, we have the following relations: $[c_\alpha^+] + [c_\alpha^-] = [c_\alpha]$ and $[c_\alpha']=2[c_\alpha]$. On the other hand the curve $c_D$ can be chosen to be:
\begin{enumerate}
\item one of the $c_\alpha$, if $\alpha\notin \Sigma_X$ and $2\alpha\notin\Sigma_X$ where $\alpha\in S$ moves $D$; or
\item one of the $c_\alpha^+$ or $c_\alpha^-$, if $\alpha\in\Sigma_X$ and $\alpha\in S$ moves $D$; or
\item $c_\alpha'$ if $\alpha\in S$ moves $D$ and $2\alpha\in\Sigma_X$.
\end{enumerate}
Now the proposition follows from the fact that the multiplicity of the intersection between $D$ and $Z$ is $\langle [D], [c_\alpha]\rangle$.
\end{proof}

Therefore we have:
\[ \dot{w}_0\delta|_{Z} = \sum_{D\in\Delta_X} n_D m_D \cdot \dot{w}_0(D\cap Z) \]
where $m_D =2$ if $D$ is moved by a simple root $\alpha\in S\cap \frac{1}{2}\Sigma$, and $m_D =1$ otherwise.

We have some last remarks which will help determining $\dot{w}_0\sigma_\mL$. We recall that $H$ is the stabilizer of a point in the open $B$-orbit.

\begin{lem} \label{lemma:appear}
If $H=N_GH$, then the coordinate $y$ must appear in the expression of $\dot{w}_0\sigma_\mL$. Moreover there exist a non-trivial line bundle $\mL_0$ generated by global sections such that $\dot{w}_0\sigma_{\mL_0}$ doesn't have the form $f(x,y^n)$ for an $n>1$.
\end{lem}
\begin{proof}
If $H=N_GH$ then proposition \ref{prop:brionemb} states that there exist a $G$-equivariant morphism $X\to\p(M)$ where $M$ is a simple $G$-module and the image of $X$ is the normalization of $X$ in its own field of rational functions. In particular, this morphism is birational onto its image. 

On the other hand, this map must be $F_{\mL_0}$ for some $\mL_0$ generated by global sections. Moreover, $F_{\mL_0}$ is given in coordinates by sections of the form $u \dot{w}_0 \sigma_{\mL_0}$ for $u\in R^u(P_X)$ (see proof of lemma \ref{lemma:diffcond}). The function $y\in\C[X_{Z,B}]$ viewed as a rational function on $X$ is stable under the action of $R^u(P_X)$.

This analysis implies that if $y$ does not appear in $\dot{w}_0 \sigma_{\mL_0}$ then the image of $X$ through $F_{\mL_0}$ is equal to the image of $Z$, which is a contradiction. This implies also that $y$ appears in $\dot{w}_0 \sigma_{\OO(D)}$ for some colour $D$, and thus in $\dot{w}_0 \sigma_{\mL}$ for all ample $\mL$.

Finally, if $\sigma_{\mL_0}$ had form $f(x,y^n)$ for an $n>1$, then all rational functions defining $F_{\mL_0}$ would share this property. This would be again in contradiction with the fact that $F_{\mL_0}$ is birational onto its image.
\end{proof}

At this point we have collected all general informations on $\dot{w}_0\sigma_\mL$ and we must examine our cuspidal rank $1$ varieties one by one. We use the list in \cite{Wa96}, and maintain the notations of this paper: we recall that we work on $\overline G$-varieties.

\subsection{Rank $1$ varieties with ${H\neq N_GH}$}
\label{ssect:nonrigid}
As one can see from the list in \cite{Wa96}, these cases are those which do not satisfy the condition (R) of theorem \ref{thm:SI}.

They are: $\mathbf {(1\mathsf A,n=2})$, $\mathbf{(7\mathsf B)}$, $\mathbf{(7\mathsf C, n=2)}$ (which is actually $\mathbf{(7\mathsf B, n=2)}$) and $\mathbf{(13)}$: we have to prove that there is no line bundle $\mL$ such that $F_\mL$ is a closed immersion.

We use the notations of section \ref{subsect:linebundles}. We have seen that the regular function $f_\delta$ on $G$ is a $H$-eigenvector with respect to the right translation, and that the associated morphism $F_\mL\colon X\to \p(V_\mL^*)$ is given by rational functions on $X$ of the form $g f_\delta/f_\delta$ for some element $g\in G$. Consider the $G$-morphism $G/H\to G/N_GH$; it is finite of degree $|N_GH/H|$ and it extends to a morphism $X\to\widetilde X$ where $\widetilde X$ is the wonderful embedding of $G/N_GH$.

In cases $\mathbf{(7\mathsf B)}$, $\mathbf{(7\mathsf C, n=2)}$ and $\mathbf{(13)}$ we have that $f_\delta$ is also an eigenvector for $N_GH$. This is easy to prove using the fact that $X$ and $\widetilde X$ both have only one colour, and the one of $X$ is the inverse image of the one of $\widetilde X$.

This implies that the map $F_\mL|_{G/H}\colon G/H \to \p(V_\mL^*)$ factorizes through $G/N_GH$ for any ample $\mL$, and thus $F_\mL$ cannot be a closed immersion.

It remains to examine the case $\mathbf {(1\mathsf A,n=2})$, which is the only one such that $f_\delta$ is not necessarily an eigenvector for the whole $N_GH$. We have $G=SL_2$, $\gamma=\alpha_1$, there are two colours $D^+,D^-$ and $S^p_X=\emptyset$. Here $X=\p^1\times\p^1$, and $G$ acts diagonally; $R^u(P_X)$ is one-dimensional. The closed $G$-orbit is $\p^1\subset X$ embedded diagonally, and is has only one colour $E$; the equation on $R^u(P_X)$ of $\dot{w}_0E$ is $x_1$. So the functions $y$ and $x_1$ have the same $T$-weight $\alpha_1$, and we can conclude that the equations of $\dot{w}_0D^+\cap X_{Z,B}$ and $\dot{w}_0D^-\cap X_{Z,B}$ both have degree $1$ in $y$, with constant leading coefficient.

Consider the line bundle $\mL= \OO(n_+D^+ + n_-D^-)$, with $n_+,n_->0$:
\[ \dot{w}_0 \sigma_\mL = (ay + x_1)^{n_+}(by + x_1)^{n_-} \]
with $a\neq b\in\C$. Now $R^u(P_X)$ is the additive group $\C$, so if $u\in R^u(P_X)$:
\[ u\cdot \dot{w}_0\sigma_\mL (x_1,y) = (ay + x_1 - u)^{n_+}(by + x_1 - u)^{n_-} \]
The derivatives of these functions calculated in $z=(0,0)$ are:
\[
\begin{array}{ccl}
\frac{\displaystyle\partial(u\cdot \dot{w}_0\sigma_\mL)}{\displaystyle\partial x_1} (0,0) &  =  &  2(-u)^{n_++n_--1}   \\[10pt]
\frac{\displaystyle\partial(u\cdot \dot{w}_0\sigma_\mL)}{\displaystyle\partial y_1} (0,0) &  =  &  (a+b) (-u)^{n_++n_--1}
\end{array}
\]
Here we see that we cannot find two elements $u_1,u_2\in R^u(P_X)$ such that the Jacobian of lemma \ref{lemma:diffcond} is nondegenerate at $z=(0,0)$, no matter which $n_+,n_-$ we choose.

\subsection{Rank $1$ varieties with $ {H=N_GH}$}
These cases are those that satisfy the condition (R) of the theorem, thus we must prove that $F_\mL$ is a closed immersion into $\p(V_\mL^*)$.

Those which are complete symmetric varieties can be omitted, since for them this fact is true (see \cite{DP83}); they are: $\mathbf {(1\mathsf A, n>2)}$, $\mathbf {(2)}$, $\mathbf {(4)}$, $\mathbf {(6\mathsf A)}$, $\mathbf {(8\mathsf B)}$, $\mathbf {(7\mathsf C, n>2)}$, $\mathbf {(8\mathsf C)}$, $\mathbf {(1\mathsf D)}$, $\mathbf {(6\mathsf D)}$, $\mathbf {(6)}$, $\mathbf {(12)}$ (again we use the labels and the ordering of \cite{Wa96}).

The remaining cases are $\mathbf {(9\mathsf B)}$, $\mathbf {(11)}$, $\mathbf {(9\mathsf C)}$, $\mathbf {(14)}$, $\mathbf {(15)}$.

Consider cases $\mathbf{(11)}$ and $\mathbf{(14)}$. $X$ has only one colour $D$. Examining the weight of the function $y$ and the weight of $f_0(x)$ (see formula \ref{formula:w0sigmaL}) for an ample line bundle $\mL = \OO(lD)$, one can apply lemma \ref{lemma:appear} and immediately conclude that $f_1(x)\neq 0$. We leave the easy details to the Reader. We can apply now for these cases the following:

\begin{lem}\label{lemma:notaroot}
If $\gamma$ does not belong to $\rmPhi^+\setminus\rmPhi_{S^p_X}$, and if $f_1(x)\neq 0$, then $F_\mL$ is a closed immersion.
\end{lem}
\begin{proof}
We begin with some considerations about the condition in lemma \ref{lemma:diffcond}. Since $z$ has coordinates $(0,\ldots,0)$ in our canonical chart, this condition is equivalent to the existence of $u_1, \ldots, u_d\in R^u(P_X)$ such that the matrix:
\[
\left(
\begin{array}{ccccc}
 \left.\frac{\partial(u_1\cdot f_0)}{\partial x_1}\right|_{x=0} & \left. \frac{\partial(u_2\cdot f_0)}{\partial x_1}\right|_{x=0} & \ldots &  \left.\frac{\partial(u_{d-1}\cdot f_0)}{\partial x_1}\right|_{x=0} &  \left.\frac{\partial(u_d\cdot f_0)}{\partial x_1}\right|_{x=0}\\[6pt]

\left.\frac{\partial(u_1\cdot f_0)}{\partial x_2}\right|_{x=0} & \left.\frac{\partial(u_2\cdot f_0)}{\partial x_2}\right|_{x=0} & \ldots & \left.\frac{\partial(u_{d-1}\cdot f_0)}{\partial x_2}\right|_{x=0} & \left.\frac{\partial(u_d\cdot f_0)}{\partial x_2}\right|_{x=0}\\[6pt]

\vdots & \vdots & \ddots & \vdots & \vdots \\[6pt]

\left.\frac{\partial(u_1\cdot f_0)}{\partial x_{d-1}}\right|_{x=0} & \left.\frac{\partial(u_2\cdot f_0)}{\partial x_{d-1}}\right|_{x=0} & \ldots & \left.\frac{\partial(u_{d-1}\cdot f_0)}{\partial x_{d-1}}\right|_{x=0} & \left.\frac{\partial(u_d\cdot f_0)}{\partial x_{d-1}}\right|_{x=0}\\[6pt]

\left.(u_1 \cdot f_1)\right|_{x=0} & \left.(u_2 \cdot f_1)\right|_{x=0} & \ldots & \left.(u_{d-1} \cdot f_1)\right|_{x=0} & \left.(u_d \cdot f_1)\right|_{x=0} \\
\end{array}
\right)
\]
is nondegenerate. Notice that one can always find at least some elements $u_1,\ldots,u_{d-1}$ such that the upper left $(d-1)\times(d-1)$ minor is non-zero, because it corresponds to the map $F_\mL$ restricted to $Z$ and this map is a closed immersion.

The nondegeneracy of the above matrix is equivalent to the fact that the function $\left.f_1(u\cdot x)\right|_{x=0}= f_1(u)$ is not a linear combination of the functions:
\begin{equation} \label{formula:lindep}
\left.\frac{\partial(f_0(u\cdot x))}{\partial x_1}\right|_{x=0}, \ldots, \left.\frac{\partial(f_0(u\cdot x))}{\partial x_{d-1}}\right|_{x=0}
\end{equation}
as functions of the variable $u\in R^u(P_X)$. Now it is quite easy to prove that the function:
\[
\left.\frac{\partial(f_0(u\cdot x))}{\partial x_i}\right|_{x=0}
\]
is a $T$-eigenvector, with weight the weight of $f_0$ minus the weight of the function $x_i$.

On the other hand, the difference between the weights of $f_0$ and $f_1$ is exactly $-\gamma$.

We switch to the appropriate set of positive roots as indexes for our variables; if we have:
\[
\sum_{\rmPhi^+\setminus\rmPhi_{S^p_X}} \mu_\alpha \left.\frac{\partial(f_0(u\cdot x))}{\partial x_\alpha}\right|_{x=0} = f_1(u)
\]
then $\mu_\alpha$ is zero except for $\alpha = \gamma$. The hypothesis of the lemma is just that there is no such $\alpha$.
\end{proof}

We remain with cases $\mathbf {(9\mathsf B)}$, $\mathbf {(9\mathsf C)}$, $\mathbf {(15)}$. All of them have two colours, moved by two different simple roots; at the beginning we can treat them in a unified way.

Let $D_1$ and $D_2$ be the two colours: they intersect $Z$ with multiplicity $1$ (proposition \ref{prop:mult}); define $E_i = D_i \cap Z$ ($i=1,2$). Call $\varphi_i(x)$ the equation on $R^u(P_X)$ of $\dot{w}_0 E_i$. We have:
\[
\begin{array}{l}
\dot{w}_0 \sigma_{\OO(D_1)}\left(x,y\right)  =  \ldots + a(x) y + \varphi_1(x) \\[5pt]
\dot{w}_0 \sigma_{\OO(D_2)}\left(x,y\right)  =  \ldots + b(x) y + \varphi_2(x)
\end{array}
\]
We take $\mL=\OO(lD_1 + sD_2)$ with $l,s>0$. Using notations of formula \ref{formula:w0sigmaL} we have:
\[
\begin{array}{l}
f_0(x)  =  \varphi_1(x)^l \varphi_2(x)^s \\[5pt]
f_1(x)  =  l a(x) \varphi_1(x)^{l-1} \varphi_2(x)^s + b(x) \varphi_1(x)^l \varphi_2(x)^{s-1}
\end{array}
\]
We can repeat the considerations in the proof of lemma \ref{lemma:notaroot}: the map $F_\mL$ is a closed immersion if and only if there exists no $\mu\in\C$ such that:
\[
\mu \left.\frac{\partial(f_0(u\cdot x))}{\partial x_\gamma}\right|_{x=0} = f_1(u)
\]
where $\gamma$ is the spherical root of $X$. Using the expression of $f_1$ and $f_0$ and dividing by $\varphi_1^{l-1}\varphi_2^{s-1}$, this equation becomes:
\[
\mu\left(l \varphi_2(u) \left.\frac{\partial(\varphi_1(u\cdot x))}{\partial x_\gamma}\right|_{x=0} + s \varphi_1(u) \left.\frac{\partial(\varphi_2(u\cdot x))}{\partial x_\gamma}\right|_{x=0}  \right) = l a(u) \varphi_2(u) + s b(u) \varphi_1(u)
\]
or equivalently:
\begin{equation}\label{formula:final}
l \varphi_2(u) \left( \mu \left.\frac{\partial(\varphi_1(u\cdot x))}{\partial x_\gamma}\right|_{x=0} - a(u) \right) = s \varphi_1(u) \left(b(u)- \mu \left.\frac{\partial(\varphi_2(u\cdot x))}{\partial x_\gamma}\right|_{x=0} \right)
\end{equation}

We examine our three cases one by one and prove that the equation above is impossible.

\vspace{10pt}\noindent
{\bf CASE $\mathbf{(9\mathsf B)}$}

\noindent
Our $\overline G$ is $PSO_{2n+1}(\C)$ ($n\geq 2$). The stabilizer $H$ of a point in the open $G$-orbit has a Levi factor isomorphic to $GL_n$, and its unipotent radical has Lie algebra isomorphic to $\bigwedge^2\C^n$ as a $GL_n$-module.

If $e_1, \ldots, e_{2n+1}$ is the canonical basis of $\C^{2n+1}$, we choose the symmetric form $(\cdot,\cdot)$ that defines the group $PSO_{2n+1}$ to be given by: $(e_i, e_j)= 1$ if $j=2n+2-i$, $(e_i,e_j)=0$ otherwise. With this choice, the Lie algebra $\mathfrak{so}_{2n+1}$ is the set of matrices that are skew symmetric around the skew diagonal. We can choose $B$ to be the (classes of) upper triangular matrices in $G$, and $T$ the (classes of) diagonal matrices.

Call $\alpha_1, \ldots, \alpha_n$ the simple roots associated to $B$ and $T$, numbered in the usual way. The spherical root of $X$ is:
\[ \gamma=\alpha_1+\ldots+\alpha_n \]

The two colours $D_1, D_2$ are moved resp. by $\alpha_1$ and $\alpha_n$. The corresponding functions $\varphi_1(x)$ and $\varphi_2(x)$ are irreducible polynomials, of weights resp.:
\[
\begin{array}{rcl}
\dot{w}_0(\omega_1)-\omega_1 & = & - 2\gamma \\[5pt]
\dot{w}_0(\omega_n)-\omega_n & = & - (\alpha_1 + 2\alpha_2 + \ldots + n\alpha_n)
\end{array}
\]
This gives a more precise information on $\dot{w}_0\sigma_{\OO(D_1)}$ and $\dot{w}_0\sigma_{\OO(D_2)}$:
\[
\begin{array}{l}
\dot{w}_0 \sigma_{\OO(D_1)}\left(x,y\right)  =  c y^2 + a(x) y + \varphi_1(x) \\[5pt]
\dot{w}_0 \sigma_{\OO(D_2)}\left(x,y\right)  =  b(x) y + \varphi_2(x)
\end{array}
\]
where $c\in\C$. We see here that $a(x)$ and $b(x)$ cannot be both zero due to lemma \ref{lemma:appear}.

The weight of $a(x)$ is $\dot{w}_0(\omega_1)-\omega_1 + \gamma$, and the weight of $b(x)$ is $\dot{w}_0(\omega_n)-\omega_n + \gamma$. This, and the irreducibiliy of $\varphi_1(x)$ and $\varphi_2(x)$, imply that if equation \ref{formula:final} is true then both sides are zero.

We examine the functions $a(x)$ and $\varphi_1(x)$. Here $P_X$ is $G_{S\setminus\{\alpha_1,\alpha_n\}}$, but the colour $D_1$ is stable under the bigger parabolic subgroup $G_{S\setminus \{\alpha_1\}}$. In the same way $\dot{w}_0 D_1$ is stable under $G_{- S\setminus \{\alpha_1\}}$. This means that $\dot{w}_0 \sigma_{\OO(D_1)}\left(x,y\right)$ will be invariant under the intersection of $R^u(P_X)$ with these two parabolic subgroups, and the function $a(x)$ too.

Denote with $K$ this intersection. If we represent $\mathsf x\in\Lie(R^u(P_X))$ as a $(2n+1)\times(2n+1)$ upper triangular matrix, then the coordinates on the first row are:
\begin{footnotesize}
\[
\begin{array} {lcl}
x_1 & = & x_{\alpha_1} \\
x_2 & = & x_{\alpha_1 + \alpha_2} \\
        & \vdots & \\
x_n & = & x_{\alpha_1+\alpha_2+\ldots+\alpha_{n-1}+\alpha_n} = x_\gamma
\end{array}
\;\;\;
\begin{array}{lcl}
x_{n+1} & = & x_{\alpha_1+\alpha_2+\ldots+\alpha_{n-1}+2\alpha_n} \\
x_{n+2} & = & x_{\alpha_1+\alpha_2+\ldots+2\alpha_{n-1}+2\alpha_n} \\
            & \vdots & \\
x_{2n-1}& = & x_{\alpha_1+2\alpha_2+\ldots+2\alpha_{n-1}+2\alpha_n}
\end{array}
\]
\end{footnotesize}
Label with $x_{2n}$, $x_{2n+1}$, etc. the remaining coordinates. With these notations we have:
\[ \Lie(K) = \left\{ x_i = 0 \;|\; i =1,\ldots\ 2n-1 \right\} \]

The weight of $a(x)$ and its invariance under $K$ tell us that:
\[ a(x) = \nu x_\gamma \]
for $\nu\in \C$.

In order to find the function $\varphi_1(x)$, we use the analysis contained in the proof of lemma \ref{lemma:weightsonZ}. The function $f_{D_1\cap Z}$ (with the notation of that lemma) is a regular function on $PSO_{2n+1}$ which is a $B$-eigenvector under left translation and a $B_-$-eigenvector under right translation; moreover $D_1$ is stable under the action of $G_{S\setminus \{\alpha_1\}}$. Thanks to these facts we find easily that $\varphi_1(x)$ is the matrix entry in the upper right corner, if we represent $x\in R^u(P_X)$ as the class of a upper triangular matrix in $SO_{2n+1}$ having all $1$'s on the diagonal. The expression in our coordinates is more complicated:
\[ \varphi_1(x) = 2 x_1x_{2n-1} + 2 x_2x_{2n-2} + \ldots + x_n^2 + \mathrm{ terms\; of\; higher\; degree }\]
The exponential $\Lie(R^u(P_X))\to R^u(P_X)$ can be expressed explicitly using the fact that matrices in $\Lie(R^u(P_X))$ are nilpotent; we leave this exercise to the Reader. Now we can write the multiplication on $R^u(P_X)$ in terms of our coordinates, and it is not difficult to conclude that:
\[ \left.\frac{\partial(\varphi_1(u\cdot x))}{\partial x_\gamma}\right|_{x=0} = u_{\alpha_1}u_{\gamma-\alpha_1} + u_{\alpha_1+\alpha_2} u_{\gamma-\alpha_1-\alpha_2}+\ldots+u_{\gamma-\alpha_n}u_{\alpha_n}  + 2 u_\gamma \]
The above expression is not a scalar multiple of $a(u)$. We conclude that the left hand side of equation \ref{formula:final} is zero if and only if $\mu= \nu =0$. But now $b$ cannot be zero, and hence the right hand side of equation \ref{formula:final} cannot be zero.

\vspace{10pt}\noindent
{\bf CASE $\mathbf{(9\mathsf C)}$}

\noindent
This case is quite similar to the one above. We have $\overline G= PSp_{2m}$, $m\geq 3$ (the case $m=2$ being equal to case $\mathbf{(9\mathsf B)}$). The stabilizer $H$ of a point in the open $G$-orbit has a Levi factor isomorphic to $\C^\times\times Sp_{2n-2}$ (up to a central isogeny), and its unipotent radical has Lie algebra isomorphic to $\C$ as a $\C^\times$-module, where $Sp_{2n-2}$ acts trivially.

If $e_1, \ldots, e_{2n}$ is the canonical basis of $\C^{2n}$, we choose the skew symmetric form $(\cdot,\cdot)$ which defines the group $PSp_{2n}$ to be given by: $(e_i, e_j)=1$ for  $j=2n+1-i$ and $1\leq i \leq n$, $(e_i, e_j)=-1$ for  $j=2n+1-i$ and $n+1\leq i \leq 2n$, $(e_i,e_j)=0$ otherwise. With this choice, the Lie algebra $\mathfrak{sp}_{2n}$ has the following form:
\[
\mathfrak{sp}_{2n}=\left\{\left(\begin{array}{cc}A & B \\C & \widetilde A\end{array}\right)\right\}\]
where $A$ is any $n\times n$-matrix, $B$ and $C$ are $n\times n$-matrices symmetric around the skew diagonal, and $\widetilde A$ is the transpose of $-A$ around the skew diagonal. We can choose $B$ to be the (classes of) upper triangular matrices in $\overline G$, and $T$ the (classes of) diagonal matrices.

Call $\alpha_1, \ldots, \alpha_n$ the simple roots associated to $B$ and $T$, numbered in the usual way. The spherical root is:
\[ \gamma=\alpha_1 + 2\alpha_2 + \ldots + 2\alpha_{n-1} + \alpha_n \]

The colours $D_1$, $D_2$ are moved resp. by $\alpha_1$ and $\alpha_2$; the corresponding functions $\varphi_1(x)$ and $\varphi_2(x)$ have weights resp.:
\[
\begin{array}{rcl}
\dot{w}_0(\omega_1)-\omega_1 & = & - \gamma - \alpha_1 \\[5pt]
\dot{w}_0(\omega_2)-\omega_2 & = & - 2\gamma
\end{array}
\]

Again, these weights provide a more precise expression of the following functions:
\[
\begin{array}{rcl}
\dot{w}_0 \sigma_{\OO(D_1)}\left(x,y\right) & = & a(x) y + \varphi_1(x) \\[5pt]
\dot{w}_0 \sigma_{\OO(D_2)}\left(x,y\right) & = & c y^2 + b(x) y + \varphi_2(x)
\end{array}
\]
where $c\in \C$, and $a(x)$, $b(x)$ are not both zero. As in the previous case, if equation \ref{formula:final} is true then both sides are zero.

The coordinates on $R^u(P_X)$ are:
\begin{footnotesize}
\[
\begin{array}{lcl}
x_1 &=& x_{\alpha_1} \\
x_2 &=& x_{\alpha_1+\alpha_2} \\
    &\vdots&  \\
x_n &=& x_{\alpha_1+\ldots+\alpha_n} \\
x_{n+1} &=& x_{\alpha_1+\ldots+\alpha_{n-2}+2\alpha_{n-1}+\alpha_n} \\
        &\vdots& \\
x_{2n-2} &=& x_{\gamma} \\
x_{2n-1} &=& x_{\gamma+\alpha_1}
\end{array}
\;\;\;
\begin{array}{lcl}
x_{2n} &=& x_{\alpha_2} \\
x_{2n+1} &=& x_{\alpha_2+\alpha_3} \\
    &\vdots&  \\
x_{3n-2} &=& x_{\alpha_2+\ldots+\alpha_n} \\
x_{3n-1} &=& x_{\alpha_2+\ldots+\alpha_{n-2}+2\alpha_{n-1}+\alpha_n} \\
        &\vdots& \\
x_{4n-4} &=& x_{\gamma-\alpha_1} \\
\end{array}
\]
\end{footnotesize}

Using the same technique as in the previous case, we find that $b$ must be invariant under the subgroup of $R^u(P_X)$ given by $x_i = 0 \; \forall i\neq 1$. This implies that:
\[ b(x) = \nu (x_\gamma - x_{\alpha_1} x_{\gamma-\alpha_1}) + \widetilde b(x) \]
where $\nu\in\C$ and $\widetilde b(x)$ is a polynomial that does not depend on the coordinates $u_\gamma$ and $u_{\gamma-\alpha_1}$. The function $\varphi_2(x)$ is the upper right $2\times2$-minor of $x\in R^u(P_X)$, and in our coordinates we have:
\[ \left.\frac{\partial(\varphi_2(u\cdot x))}{\partial x_\gamma}\right|_{x=0} = 2u_\gamma - u_{\alpha_1} u_{\gamma-\alpha_1} \]
The above expression is not a scalar multiple of $b(u)$, hence $b$ and $\mu$ are zero. But now $a$ cannot be zero, nor the left hand side of equation \ref{formula:final}.

\vspace{10pt}\noindent
{\bf CASE $\mathbf{(15)}$}

\noindent
Here our group $G$ is of type $\mathsf G_2$, $H$ has a Levi factor quotient of $\C^\times \times SL_2$, and the Lie algebra of its unipotent radical is isomorphic to $\C \oplus \C^2$ as a $\C^\times \times SL_2$-module. If $\alpha_1$ and $\alpha_2$ are the two simple roots (short and long, resp.), the spherical root of $X$ is:
\[ \gamma = \alpha_1 + \alpha_2 \]

The colours $D_1$ and $D_2$ are moved resp. by $\alpha_1$ and $\alpha_2$, so $B=P_X$. The functions $\varphi_1$ and $\varphi_2$ have weights:
\[
\begin{array}{rcl}
\dot{w}_0(\omega_1)-\omega_1 & = & - 4\alpha_1 - 2\alpha_2 \\[5pt]
\dot{w}_0(\omega_2)-\omega_2 & = & - 6\alpha_1 - 4\alpha_2
\end{array}
\]
hence:
\[
\begin{array}{rcl}
\dot{w}_0 \sigma_{\OO(D_1)}\left(x,y\right) & = & c(x) y^2 + a(x) y + \varphi_1(x) \\[5pt]
\dot{w}_0 \sigma_{\OO(D_2)}\left(x,y\right) & = & f(x) y^4 + e(x) y^3 + d(x) y^2 + b(x) y + \varphi_2(x)
\end{array}
\]

In order to do explicit calculations on $R^u(P_X)$, we use the embedding of $\mathsf G_2$ inside $SO_8(\C)$ given by:
\begin{footnotesize}
\[ \Lie(G) = \left\{ \left(
\begin{array}{cccccccc}
a_1    & x_7     & x_2     & -x_4 & x_4  & x_5     & x_3      & 0     \\
x_6    & a_1+a_2 & x_4     & -x_5 & x_5  & x_1     & 0        & -x_3  \\
x_{11} & x_9     & a_2     & x_6  & -x_6 & 0       & -x_1     & -x_5  \\
-x_9   & -x_8    & x_7     & 0    & 0    & x_6     & -x_5     & -x_4  \\
x_9    & x_8     & -x_7    & 0    & 0    & -x_6    & x_5      & x_4   \\
x_8    & x_{12}  & 0       & x_7  & -x_7 & -a_2    & -x_4     & -x_2  \\
x_{10} & 0       & -x_{12} & -x_8 & x_8  & -x_9    & -a_1-a_2 & -x_7  \\
0      & -x_{10} & -x_8    & -x_9 & x_9  & -x_{11} & -x_6     & -a_1
\end{array}
\right)\right\} \]
\end{footnotesize}
for $a_1$, $a_2$, $x_1$, \ldots, $x_{12}$ coordinates on $\Lie(G)$. We have:
\begin{footnotesize}
\[ \Lie(R^u(Q)) = \left\{ \left(
\begin{array}{cccccccc}
0   & 0   & x_2 & -x_4 & x_4  & x_5  & x_3  & 0     \\
x_6 & 0   & x_4 & -x_5 & x_5  & x_1  & 0    & -x_3  \\
0   & 0   & 0   & x_6  & -x_6 & 0    & -x_1 & -x_5  \\
0   & 0   & 0   & 0    & 0    & x_6  & -x_5 & -x_4  \\
0   & 0   & 0   & 0    & 0    & -x_6 & x_5  & x_4   \\
0   & 0   & 0   & 0    & 0    & 0    & -x_4 & -x_2  \\
0   & 0   & 0   & 0    & 0    & 0    & 0    & 0     \\
0   & 0   & 0   & 0    & 0    & 0    & -x_6 & 0     \\
\end{array}
\right)\right\} \]
\end{footnotesize}
These $x_i$ are our coordinates given by root space decomposition, and if we label them using the corresponding roots, we have:
\[
\begin{array}{lcl}
x_1 &=& x_{3\alpha_1+\alpha_2} \\
x_2 &=& x_{\alpha_2} \\
x_3 &=& x_{3\alpha_1+2\alpha_2}
\end{array}
\;\;\;
\begin{array}{lcl}
x_4 &=& x_{\alpha_1+\alpha_2} \\
x_5 &=& x_{2\alpha_1+\alpha_2} \\
x_6 &=& x_{\alpha_1}
\end{array}
\]
At this point we need the expression of the group operation on $R^u(P_X)$ in terms of these coordinates. This is an easy task since $\Lie(R^u(P_X))$ is nilpotent and thus the exponential map can be written explicitly; we leave this exercise to the Reader.

As in case $\mathbf{(9\mathsf B)}$, the functions $a(x)$, $b(x)$, $c(x)$, $d(x)$, $e(x)$ and $f(x)$ are stable under some nontrivial subgroup of $R^u(P_X)$.

The possible monomials occurring in $a(x)$ are $x_1$, $x_5 x_6$, $x_4 x_6^2$, $x_2 x_6^3$. Invariance under the translation by $\exp(\mathfrak g_{\alpha_2})$ tell us that:
\[ a(x) = \mu_1 x_1 + (3\mu_2 - 6\mu_3) x_5 x_6 + \mu_2 x_4 x_6^2 + \mu_3 x_2 x_6^3 \]
with $\mu_1, \mu_2, \mu_3\in \C$.

There exist sixteen possible monomials occurring in $b(x)$. Unfortunately, lemma \ref{lemma:appear} does not help proving that $a(x)$ or $b(x)$ are non-zero, due to the presence of $e(x)$ (which one can prove it is actually non-zero!). We need a deeper analysis, which we describe leaving the computations to the Reader.

We set up a system of coordinates $\xi_1,\ldots,\xi_{14}$ on the big cell $R^u(B) T R^u(B_-)\cong \C^6\times \left(\C^\times\right)^2 \times \C^6$ of $G$, using the exponential map for $R^u(B)$ and $R^u(B_-)$. We find the subgroup $H$ described in \cite{Wa96}, and we can express it in terms of the coordinates $\xi_1,\ldots,\xi_{14}$.

In \cite{Wa96} we find the $B$-weight and the $H$-weight\footnote{As usual, $B$ acts on functions on $G$ by left translation and $H$ acts by right translation.} of $f_{D_2}\in\C[G]$ a global equation of $D_2$ pulled back to $G$ along the projection $\pi\colon G\to G/H$. These weights are enough to find $f_{D_2}$ up to a multiplicative constant, using the decomposition of $\C[G]$ as a $G\times G$-module. We obtain the rational function $F=\dot{w}_0 f_{D_2}/f_{D_2}\in\C(G/H)\subset \C(G)$. The function $F$ is nothing but $\dot{w}_0 \sigma_{\OO(D_2)}$ expressed in our coordinates on the big cell of $G$.

Fix an arbitrary point $p$ in the big cell, lying also inside $BH$. Since $p\in BH$, then $\pi(p)$ is inside the canonical chart $X_{Z,B}$ of $X$, as well as the whole $\pi(u \cdot p)$ for all $u\in R^u(P_X)$. More precisely, if $y_0\in\C$ is the value of the coordinate $y$ in the point $\pi(p)$, then:
\[\left\{ \pi(u \cdot p) \;\;|\;\; u\in R^u(Q)\right\} = \{y = y_0 \}\subset X_{Z,B} \]

Consider the function $F(u \cdot p)$ for $u\in R^u(P_X)$, expressed in terms of the coordinates $x(u)=(x_1(u),\ldots,x_{d-1}(u))$ of $u$. If $\pi(p)$ had all coordinates $x_1,\ldots,x_{d-1}$ equal to zero, then $F(u\cdot p)$ would actually be $\dot{w}_0 \sigma_{\OO(D_2)}(x(u),y_0)$, and we could recover the functions $f(x)$, $e(x)$, $d(x)$, $b(x)$ and $\varphi_2(x)$ from it. The problem is that we have no control on the coordinates of $\pi(p)$ relative to $X_{Z,B}$, so we cannot go back and choose $p$ in  order to make this happen.

However, since $R^u(P_X)$ is unipotent, we have:
\[ F(u \cdot p) = \dot{w}_0 \sigma_{\OO(D_2)}(x(u),y_0) + \mathrm{ other \;terms\; depending\; on\; the\; coordinates\; of\; } p \]

Now, the explicit calculations show that the possible monomials of $b(x)$ appear in $F(u \cdot p)$ with coefficients which depend only on $y_0$ (and not on the other coordinates of $p$). This implies that from the expression of $F(u\cdot p)$ we can actually recover $b(x)$:
\[
\begin{array}{rl}
b(x) =& \mu_4 \left(-720 x_3x_5 + 720 x_3x_4x_6 -240 x_3x_6^2x_2 + 360 x_1x_5x_2 -720 x_1x_4^2 \right.\\[5pt]
      & 360 x_1x_4x_6x_2  -60 x_1x_6^2x_2^2 + 360 x_5^2x_4  -360  x_5^2x_6x_2 + \\[5pt]
      & \left. 60 x_5x_4x_6^2x_2 + 18 x_5x_6^3x_2^2 -8 x_4x_6^4x_2^2 + x_6^5x_2^3\right)
\end{array}
\]
for $\mu_4\in \C\setminus \{0\}$. Therefore $b\neq 0$.

Functions $\varphi_1$ and $\varphi_2$ can be found as in the previous cases, thanks to the inclusion $\mathsf G_2\subset SO_8$. Indeed, we can use the functions on $SO_8$ which are $\widetilde B$-eigenvectors under left translation and $\widetilde B_-$-eigenvectors under right translation, where $\widetilde B$ is a suitable Borel subgroup of $SO_8$. For $x\in R^u(P_X)\subset SO_8$, the function $\varphi_1(x)$ is the pfaffian of the upper right $3\times3$-submatrix of $x$, and $\varphi_2(x)$ is the upper right $2\times2$-minor of $x$. The expressions in our coordinates are:
\[
\begin{array}{rl}
\varphi_1(x)=& 360 x_1 x_4 + 360 x_5^2 - 360 x_3 x_6 + 30 x_2 x_5 x_6^2 - x_2^2 x_6^4 \\[5pt]
\varphi_2(x)=& \frac{1}{4} x_1^2x_2^2 -x_3^2 - \frac{3}{4} x_4^2 x_5^2 +x_2x_5^3 +x_1\left(x_4^3 - \frac{3}{2}x_2x_4x_5\right) + \\[5pt]
             & \frac{1}{240} \left( x_2^2 x_4^2 x_6^4 - x_2^3 x_5 x_6^4\right) - \frac{1}{21600}x_2^4x_6^6 + x_3\left(-x_4^2x_6 + x_2x_5x_6 + \frac{1}{10}x_2^2x_6^3 \right)
\end{array}
\]

At this point, the formulas we have found and the multiplication on $R^u(P_X)$ in our coordinates are enough to express explicitly equation \ref{formula:final}; the proof that this equation cannot be true is straightforward.

The results in this section, together with lemma \ref{lemma:reduction}, complete the proof of theorem \ref{thm:SI}.

\section{Proof of Theorem \ref{thm:main}}
\label{sect:proof}
We begin proving the uniqueness of a $G$-equivariant closed immersion $F\colon X\to \p(V)$ for any fixed simple $V$, using theorem \ref{thm:SI}.

\begin{lem}
\label{lemma:restriction}
Let $X$ be a wonderful variety admitting a $G$-equivariant closed immersion into the projective space of a simple $G$-module. Then the restriction of a line bundle from $X$ to $Z$, the unique closed orbit, gives an inclusion $\Pic(X)\subseteq \Pic(Z)$ where the latter is identified as usual with a sublattice of the integral weights.
\end{lem}
\begin{proof}
Theorem \ref{thm:SI} guarantees that inside $X$ there is no $G$-stable subvariety being a parabolic induction of the $SL_2$-variety $\p^1\times\p^1$: this implies that there is no simple root moving two colours on $X$. By \cite{Lu01} (proposition 3.2), it follows that any colour $D$ on $X$ will be moved either by a single simple root $\alpha_D$, or by two orthogonal simple roots $\alpha_D, \alpha'_D$.

In view of proposition \ref{prop:mult}, a colour $D$ will intersect $Z$ in the union of at most two colours of $Z$, in such a way that $D\cap Z$ corresponds to either $\omega_D$, or $2\omega_D$, or $\omega_D+\omega'_D$ (where $\omega_D$ is the fundamental dominant weight corresponding to $\alpha_D$). These three cases occur resp. when $D$ is moved by $\alpha_D\in S$ with $2\alpha_D \notin \Sigma_X$, or $D$ is moved by $\alpha_D\in S$ with $2\alpha_D\in\Sigma_X$, or $D$ is moved by $\alpha_D,\alpha'_D \in S$. In this situation, moreover, two different colours will be moved by two disjoint sets of simple roots, and this completes the proof. 
\end{proof}

The identification of $\Pic(X)$ with a sublattice of the weights is exactly the application $\mL \mapsto \chi_\mL$. Fix $V$ and suppose we are given a closed immersion $F\colon X\to \p(V)$: thanks to the lemma, the highest weight of $V$ determines uniquely the line bundle giving this immersion. More precisely, there exists a unique line bundle $\mL$ such that $V\cong V_{\mL}^*$. This implies that $F$ is determined up to composition with elements of $GL(V)$. But we are dealing with $G$-equivariant maps: any $A\in GL(V)$ such that $A\circ F$ remains $G$-equivariant will have to commute with the $G$-action on $V$, thus $A$ will act as a scalar on $V$. This proves uniqueness of $F$.

\vspace{10pt}
\noindent
{\bf Remark}
The image of $\mL \mapsto \chi_\mL$ (for $\mL$ varying among the ample line bundles on $X$) gives a subset of dominant weights which determines exactly all the simple $G$-modules whose projective space contains a copy of $X$.
\vspace{10pt}

Finally, what we are left to prove is that if $X$ has property (R) of theorem \ref{thm:SI} then it is strict, i.e. all the stabilizers of its points are equal to their normalizers. Let $x\in X$ and let $G_x$ be its stabilizer. The closure of the $G$-orbit of $x$ is a wonderful $G$-subvariety $Y$, whose generic stabilizer is $G_x$.

This $Y$ has the property (R) as well, so there exists a simple $G$-module $V$ and a unique closed immersion $F\colon Y\to\p(V)$. Suppose that $G_x$ is different from its normalizer, and take an element in $N_G(G_x)\setminus G_x$. This element induces a non-trivial $G$-equivariant automorphism $\phi$ of $Y$; this is absurd because then $F$ and $F\circ\phi$ would be two different closed immersions of $Y$ in $\p(V)$. This finishes the proof of theorem \ref{thm:main}.

\section{Ample and very ample line bundles}
\label{sect:veryample}
\begin{thm}
\label{thm:veryample}
Any ample line bundle on a wonderful variety is very ample.
\end{thm}
\begin{proof}
For strict wonderful varieties, theorem \ref{thm:SI} assures our result. For non-strict varieties, the problem can be reduced to rank $1$ cuspidal wonderful varieties exactly in the same way as in lemma \ref{lemma:reduction}. We remark that here we cannot ignore the cases where the centre of $G$ doesn't act trivially.

Following again \cite{Wa96}, the non-strict cuspidal wonderful varieties of rank $1$ are: $\mathbf{(1\mathsf A)}$,  $\mathbf{(3)}$, $\mathbf{(5\mathsf A)}$, $\mathbf{(7\mathsf B)}$, $\mathbf{(10)}$, $\mathbf{(5\mathsf D)}$, $\mathbf{(5)}$, $\mathbf{(13)}$. Some of these varieties have easy explicit descriptions, however there is no need of a case-by-case proof. Let $X$ be one of these varieties, $\gamma$ the spherical root of $X$, and $f_\gamma$ a rational function on $X$, $B$-eigenvector of weight $\gamma$. From table 1 in \cite{Wa96}, we see that in all these cases $1/f_\gamma$ has poles of order $1$ on (all) the colour(s) of $X$, a zero of order $1$ on the closed orbit (by construction), and no other poles.

Therefore $1/f_\gamma$ belongs to $\Gamma(X,\mL)$ for any ample line bundle $\mL$ associated to a sum of colours with positive coefficients. If we focus on the canonical chart as in section \ref{sect:rank1}, the function $1/f_\gamma$ is nothing but the function $y\in\C[X_{Z,B}]$ (up to a multiplicative constant).

Let us use the notations of lemma \ref{lemma:diffcond}; thanks to the Borel-Weil theorem there exist $u_1,\ldots,u_{d-1}\in R^u(P_X)$ such that the jacobian matrix of the functions $(u_1\dot{w}_0)\sigma_\mL$, $\ldots$, $(u_{d-1}\dot{w}_0)\sigma_\mL$ with respect only to the coordinates of $R^u(P_X)$ is nondegenerate in $z$. But now $y$ is among the global sections we can consider, and it is clear that the jacobian matrix of the functions $(u_1\dot{w}_0)\sigma_\mL,\ldots, (u_{d-1}\dot{w}_0)\sigma_\mL, y$ with respect to all the coordinates is nondegenerate in $z$.

Therefore in $\Gamma(X,\mL)$ there are enough sections to give an immersion $X\to\p(\Gamma(X,\mL)^*)$, and $\mL$ is very ample.
\end{proof}

\end{document}